\documentclass[12pt]{article}
\usepackage{amsmath,amssymb,amsthm, comment}

\newtheorem{theorem}{Theorem}[section]
\newtheorem{thmy}{Theorem}

\newtheorem{corollary}[theorem]{Corollary}

\def\barr{\begin{array}}
\def\earr{\end{array}}

\title{A result on certain sums of element orders in finite groups}
\author{Marius T\u arn\u auceanu}
\date{April 2, 2025}

\begin{document}

\maketitle

\begin{abstract}
Given a finite group $G$ of order $p^nm$, where $p$ is a prime and $p\nmid m$, we denote by $\psi_p(G)$ the sum of orders of $p$-parts of elements in $G$. In the current note, we prove that $\psi_p(G)\leq\psi_p(C_{p^nm})$, where $C_{p^nm}$ is the cyclic group of order $p^nm$, and the equality holds if and only if $G$ is $p$-nilpotent of a particular type. A generalization of this result is also presented.
\end{abstract}

{\small
\noindent
{\bf MSC2020\,:} Primary 20D60; Secondary 20D15, 20F18.

\noindent
{\bf Key words\,:} element orders, finite groups.}

\section{Introduction}
Let $G$ be a finite group. In 2009, H. Amiri, S.M. Jafarian Amiri and I.M. Isaacs introduced in their paper \cite{1} the function
\begin{equation}
\psi(G)=\sum_{x\in G}o(x),\nonumber
\end{equation}where $o(x)$ denotes the order of $x$ in $G$. They proved the following basic theorem:

\begin{thmy}
If $G$ is a finite group of order $n$, then $\psi(G)\leq\psi(C_n)$, and we have equality if and only if $G$ is cyclic.
\end{thmy}Since then many authors have studied the properties of the function $\psi(G)$ and its relations with the structure of $G$ (see e.g. \cite{3}-\cite{5}, \cite{8}-\cite{10} and \cite{15}).

Note that Theorem A follows from the next result:

\begin{thmy}
If $G$ is a finite group of order $n$, then there is a bijection $f:G\longrightarrow C_n$ such that $o(x)$ divides $o(f(x))$, for all $x\in G$.
\end{thmy}

This has been formulated as a question by I.M. Isaacs (see Problem 18.1 in \cite{14}) and proved for some particular groups by F. Ladisch \cite{13} and M. Amiri and S.M. Jafarian Amiri \cite{2}. A proof for arbitrary groups has been recently given by M. Amiri \cite{6}. 

Next we assume that $|G|=p^nm$, where $p$ is a prime and $p\nmid m$, and consider the function
\begin{equation}
\psi_p(G)=\sum_{x\in G}o(x_p),\nonumber
\end{equation}where $x_p$ is the $p$-part of $x$ in $G$. We remark that
\begin{equation}
\psi_p(G)=\sum_{x\in G}o(x^m)\nonumber
\end{equation}and in particular
\begin{equation}
\psi_p(C_{p^nm})=m\psi(C_{p^n}).\nonumber
\end{equation}

Our main result is stated as follows.

\begin{theorem}
If $G$ is a finite group of order $p^nm$ with $p$ prime and $p\nmid m$, then $\psi_p(G)\leq\psi_p(C_{p^nm})$, and we have equality if and only if $G\cong H\rtimes C_{p^n}$, where $H$ is a normal subgroup of order $m$ of $G$.
\end{theorem}Its proof will be given in Section 2, while in Section 3 we will present a generalization of this theorem. Also, we remark that there is no constant $c\in(0,1)$ such that if $\psi_p(G)>c\psi_p(C_{p^nm})$, then $G$ belongs to a significant class of groups, as it happens for the original function $\psi$.

For the proof of Theorem 1.1, we need two well-known theorems of Frobenius \cite{7} and  N. Iiyori and H. Yamaki \cite{11}.

\begin{thmy}
Let $G$ be a finite group whose order is divisible by a number $n$. Then the number of solutions of the equation $x^n=1$ in $G$ is a
multiple of $n$.
\end{thmy}

\begin{thmy}
Let $G$ be a finite group whose order is divisible by a number $n$. If the set of solutions of the equation $x^n=1$ in $G$ has exactly $n$
elements, then it forms a subgroup of $G$.
\end{thmy}

Most of our notation is standard and will usually not be repeated here. Elementary notions and results on groups can be found in \cite{7}.

\section{Proof of Theorem 1.1}
First of all, we observe that if $G$ is a finite group of order $n$, then Theorem B implies the existence of a partition
$(L_d(G))_{d|n}$ of $G$ such that for every divisor $d$ of $n$ we have:
\begin{itemize}
\item[{\rm a)}] $|L_d(G)|=\varphi(d)$, where $\varphi$ is the Euler totient function;
\item[{\rm b)}] $x^d=1$, $\forall\, x\in L_d(G)$.
\end{itemize}Moreover, the bijection $f$ in Theorem B maps the elements of $L_d(G)$ to the elements of order $d$ of $C_n$.
\bigskip

We are now able to prove our main result.

\bigskip\noindent{\bf Proof of Theorem 1.1.} Under the above notation, it is easy to see that 
\begin{equation}
o(x)|o(f(x)) \Rightarrow o(x_p)|o(f(x)_p), \forall\, x\in G,\nonumber 
\end{equation}which leads to
\begin{align*}
\psi_p(G)&=\sum_{x\in G}o(x_p)\\
&\leq\sum_{x\in G}o(f(x)_p)\\
&=\sum_{y\in C_{p^nm}}o(y_p)\\
&=\psi_p(C_{p^nm}).\nonumber
\end{align*}

Assume now that $\psi_p(G)=\psi_p(C_{p^nm})$ and let $0\leq i\leq n$ and $m_1|m$. Then every element $x\in L_{p^im_1}(G)$ can be written as $x=x_{pi}x_i'$, where $o(x_{pi})|p^i$ and $o(x_i')|m_1$. We get
\begin{equation}
\psi_p(G)=\sum_{i=0}^n\sum_{m_1|m}\sum_{x\in L_{p^im_1}(G)}o(x_{pi})\nonumber
\end{equation}and
\begin{equation}
\psi_p(C_{p^nm})=\sum_{i=0}^n\sum_{m_1|m}\varphi(p^im_1)p^i,\nonumber
\end{equation}implying that $o(x_{pi})=p^i$, $\forall\, x\in L_{p^im_1}(G)$. Thus 
\begin{equation}
o(x_{pn})=p^n, \forall\, x\in L_{p^nm_1}(G),\nonumber
\end{equation}and so the Sylow $p$-subgroups of $G$ are cyclic. Moreover
\begin{align*}
|\{x\in G\mid x^m=1\}|&=\sum_{m_1|m}|L_{p^0m_1}(G)|\\
&=\sum_{m_1|m}\varphi(m_1)\\
&=m\nonumber
\end{align*}and Theorem D shows that $H=\{x\in G\mid x^m=1\}$ is a subgroup of $G$. Clearly, $H$ is normal in $G$ and it follows that $G\cong H\rtimes C_{p^n}$, completing the proof.\qed

We remark that an alternative way to prove the cyclicity of the Sylow $p$-subgroups of $G$ is as follows. Since 
\begin{equation}
\psi_p(G)=\psi_p(C_{p^nm})=m\psi(C_{p^n})=m\,\frac{p^{2n+1}+1}{p+1}\,,\nonumber 
\end{equation}there is $x\in G$ such that 
\begin{equation}
o(x_p)\geq\frac{p^{2n+1}+1}{p^n(p+1)}>p^{n-1},\nonumber 
\end{equation}i.e. $o(x_p)=p^n$.
\bigskip

We end this section with the following example.

\bigskip\noindent{\bf Example.} The group $G={\rm SL}_2(\mathbb{F}_3)=Q_8\rtimes C_9=SmallGroup(72,3)$ has $1$ element of order $1$, $1$ element of order $2$, $2$ elements of order $3$, $6$ elements of order $4$, $2$ elements of order $6$, $24$ elements of order $9$, $12$ elements of order $12$ and $24$ elements of order $18$. The sets $(L_d(G))_{d|12}$ can be chosen as follows
\begin{itemize}
\item[-] $L_1(G)$: $1$ element of order $1$,
\item[-] $L_2(G)$: $1$ element of order $2$,
\item[-] $L_3(G)$: $2$ elements of order $3$,
\item[-] $L_4(G)$: $2$ elements of order $4$,
\item[-] $L_6(G)$: $2$ elements of order $6$,
\item[-] $L_8(G)$: $4$ elements of order $4$,
\item[-] $L_9(G)$: $6$ elements of order $9$,
\item[-] $L_{12}(G)$: $4$ elements of order $12$,
\item[-] $L_{18}(G)$: $6$ elements of order $9$,
\item[-] $L_{24}(G)$: $8$ elements of order $12$,
\item[-] $L_{36}(G)$: $12$ elements of order $9$,
\item[-] $L_{72}(G)$: $24$ elements of order $18$,
\end{itemize}and we easily get
$$\psi_2(G)=387=\psi_2(C_{72}) \mbox{ and } \psi_3(G)=488=\psi_3(C_{72}).$$

\section{A generalization}

Let $G$ be a finite group of order $n$ and $\pi$ be a set of primes contained in $\pi(n)$. Then every element $x$ of $G$ can be uniquely written as $x=x_{\pi}x_{\pi'}$, where $o(x_{\pi})$ is a $\pi$-number and $o(x_{\pi'})$ is a $\pi'$-number. Define
\begin{equation}
\psi_{\pi}(G)=\sum_{x\in G}o(x_{\pi}).\nonumber
\end{equation}

Similarly with Theorem 1.1, one can prove the next result.

\begin{theorem}
Under the above notation, we have $\psi_{\pi}(G)\leq\psi_{\pi}(C_n)$, and the equality holds if and only if $G\cong H\rtimes C$, where $H$ is a normal subgroup of $G$ of order the $\pi'$-part of $n$ and $C$ is a cyclic subgroup of $G$ of order the $\pi$-part of $n$.
\end{theorem}

Note that an immediate consequence of Theorem 3.1 is the following.

\begin{corollary}
Under the above notation, if $\psi_{\pi}(G)=\psi_{\pi}(C_n)$ and $\psi_{\pi'}(G)=\psi_{\pi'}(C_n)$, then $G\cong C_n$.
\end{corollary}

Finally, we formulate the following question concerning Theorem D in Section 1.
\bigskip

\noindent{\bf Question.} Let $G$ be a finite group of order $n$, and let $P\in Syl_p(G)$ be a cyclic group. If $|\{x\in G : x^{n_{p'}}=1\}|=n_{p'}$ , then $\{x\in G : x^{n_{p'}}=1\}$ is a subgroup of $G$.
\bigskip

\noindent{\bf Acknowledgements.} The author is grateful to the reviewer for remarks which improve the previous version of the paper.
\bigskip

\noindent{\bf Funding.} The author did not receive support from any organization for the submitted work.
\bigskip

\noindent{\bf Conflicts of interests.} The author declares that he has no conflict of interest.
\bigskip

\noindent{\bf Data availability statement.} My manuscript has no associated data.

\vspace*{3ex}\small

\hfill
\begin{minipage}[t]{5cm}
Marius T\u arn\u auceanu \\
Faculty of  Mathematics \\
``Al.I. Cuza'' University \\
Ia\c si, Romania \\
e-mail: {\tt tarnauc@uaic.ro}
\end{minipage}

\end{document}